\documentclass{amsart}

\usepackage{amscd,amssymb,amsmath,graphicx}
\usepackage{mathrsfs}
\usepackage[all]{xy}

\newcommand\mylabel[1]{\label{#1}}

\newtheorem{theorem}{Theorem}[section]

\newtheorem{proposition}[theorem]{Proposition}

\newtheorem{conjecture}[theorem]{Conjecture}

\theoremstyle{definition}
\newtheorem{definition}[theorem]{Definition}

\theoremstyle{remark}

\DeclareFontFamily{U}{wncy}{}
\DeclareFontShape{U}{wncy}{m}{n}{<->wncyr10}{}
\DeclareSymbolFont{mcy}{U}{wncy}{m}{n}
\DeclareMathSymbol{\Sh}{\mathord}{mcy}{"58}

\newcommand{\ZZ}	{\mathbb{Z}}
\newcommand{\QQ}	{\mathbb{Q}}

\newcommand{\CC}	{\mathbb{C}}

\newcommand{\PP}	{\mathbb{P}}
\renewcommand{\AA}	{\mathbb{A}}
\newcommand{\GG}	{\mathbb{G}}

\newcommand  {\shA}     {\mathscr{A}}

\newcommand  {\shE}     {\mathscr{E}}
\newcommand  {\shF}     {\mathscr{F}}

\newcommand  {\shI}     {\mathscr{I}}

\newcommand  {\shN}     {\mathscr{N}}
\newcommand  {\shL}     {\mathscr{L}}

\newcommand  {\an}      {{\text{\rm an}}}

\newcommand  {\Aut}     {\operatorname{Aut}}

\newcommand  {\Br}      {\operatorname{Br}}

\newcommand  {\Cl}      {\operatorname{Cl}}

\newcommand  {\GL}      {\operatorname{GL}}
\newcommand  {\GO}      {\operatorname{GO}}

\newcommand  {\Hom}     {\operatorname{Hom}}
\newcommand  {\Hilb}    {\operatorname{Hilb}}

\newcommand  {\id}      {{\operatorname{id}}}

\newcommand  {\lra}     {\longrightarrow}

\renewcommand{\O}       {\mathscr{O}}

\newcommand  {\Pic}     {\operatorname{Pic}}
\newcommand  {\PGL}     {\operatorname{PGL}}

\newcommand  {\pr}      {\operatorname{pr}}
\newcommand  {\Proj}    {\operatorname{Proj}}

\newcommand  {\ra}      {\rightarrow}

\newcommand  {\res}     {\operatorname{res}}

\newcommand  {\Sing}    {\operatorname{Sing}}
\newcommand  {\SL}      {\operatorname{SL}}

\newcommand  {\SO}      {\operatorname{SO}}

\newcommand  {\Spec}    {\operatorname{Spec}}
\newcommand  {\Spin}    {\operatorname{Spin}}

\newcommand  {\Sym}     {\operatorname{Sym}}

\def\mydate{\number\day\space\ifcase\month \or January\or February\or March\or 
April\or May\or June\or July\or
August\or September\or October\or November\or December\fi \space\number\year}

\DeclareFontFamily{U}{wncy}{}
\DeclareFontShape{U}{wncy}{m}{n}{<->wncyr10}{}
\DeclareSymbolFont{mcy}{U}{wncy}{m}{n}
\DeclareMathSymbol{\Sh}{\mathord}{mcy}{"58}

\begin{document}

\title[Brauer groups for quiver moduli]
      {Brauer groups for quiver moduli}

\author[Markus Reineke]{Markus Reineke}
\address{Fachbereich C - Mathematik, Bergische Universit\"at Wuppertal, Gau\ss str.\ 20,
42097 Wuppertal,  Germany}
\curraddr{}
\email{mreineke@uni-wuppertal.de}

\author[Stefan Schr\"oer]{Stefan Schr\"oer}
\address{Mathematisches Institut, Heinrich-Heine-Universit\"at,
40204 D\"usseldorf, Germany}
\curraddr{}
\email{schroeer@math.uni-duesseldorf.de}



\begin{abstract} We compute the Brauer groups of several moduli spaces of stable quiver representations.\end{abstract}

\maketitle
\tableofcontents

\section*{Introduction}
In posing moduli problems for isomorphism classes of objects in a category, one usually aims at constructing a fine moduli space, that is, one which representents the given moduli functor, and in particular carries a universal, or at least tautological, family of objects. However, this is typically obstructed by numerical conditions on appropriate discrete invariants of the objects to be parametrized, a typical example being vector bundles on smooth projective curves, which only admit fine moduli spaces for coprimality of rank and degree of the bundles. In general, one thus only arrives at coarse moduli spaces, which, a priori, do not carry tautological families of objects.

To actually disprove existence of such tautological families on coarse moduli spaces constitutes a completely separate problem,  whose solution is closely tied to computing non-trivial geometric invariants of the coarse moduli spaces. Again taking moduli of vector bundles on curves as an example, such nonexistence proofs were given in \cite{Ramanan 1973},  \cite{Drezet; Narasimhan 1989}, \cite{Balaji; Biswas; Gabber; Nagaraj 2007} using computations of the canonical class, the Picard group, and the Brauer group of the moduli spaces, respectively.

Quiver moduli, which are moduli spaces of stable representations of quivers provide another instance of such moduli problems (see \cite{King 1994} for the construction and \cite{moduli} for an overview over their known geometric and topological properties). In fact, their behaviour is often expected to be close to the case of vector bundles on curves since in both cases the underlying category is of homological dimension one. Quiver moduli attracted interest in recent years in particular due to their role in motivic Donaldson-Thomas theory (see, for example, \cite{Reineke 2010},\cite{Reineke 2011}, \cite{Reineke; Weist 2013}).

The topic of the present work is the computation of the Brauer group of quiver moduli, taking up the approach of \cite{Balaji; Biswas; Gabber; Nagaraj 2007}, and the study of the non-existence of tautological bundles on them. We conjecture that the Brauer group of the moduli space parametrizing stable representations of a given dimension type ${\bf d}$ is cyclic of order the greatest common divisor of the entries of ${\bf d}$ (Conjecture \ref{main conj}).

Our main result, Theorem \ref{main theorem}, proves this conjecture under a mild codimension assumption. We prove that this assumption is always fulfilled for two of the most promiment classes of quivers, namely multiple loop and generalized Kronecker quivers, except for two special cases which can be treated by explicit coordinatization.

As a consequence, we prove the desired non-existence of tautological representations in the case of non-primitive dimension vectors, under the mentioned codimension assumption, and unconditionally for multiple loop and generalized Kronecker quivers.

We also compute the Brauer group for moduli of quadrics, which is a closely related moduli problem of linear algebra data. Here we obtain results over arbitrary ground rings including the case of characteristic $p=2$.

The paper is organized as follows:

We first collect all necessary facts on Brauer groups in Section \ref{brauer}. Since results seem to be scattered in the literature, we state and prove them systematically for the reader's convenience. Then we collect all facts related to quiver representations and their moduli spaces, unframed as well as framed, in Section \ref{quiver}; for a more thorough treatment, the reader is referred to \cite{moduli}. 

After these preparations, we can prove our main theorem in Section \ref{main}, adapting the strategy of \cite{Balaji; Biswas; Gabber; Nagaraj 2007} to the quiver setup. Here we also state the application to the non-existence of universal and tautological representations.

Concerning the crucial codimension assumption of the main theorem, we derive a sufficient, purely numerical, criterion for this to be fulfilled in Section \ref{combinatorics}. This derivation uses standard stratification techniques of the theory of quiver moduli, for example the Harder--Narasimhan stratification of \cite{Reineke 2003}.

In Section \ref{lkq}, we show that this numerical criterion is (almost) always fulfilled for multiple loop and generalized Kronecker quivers. The general case is elementary in nature, but requires a rather technical case-by-case analysis. The two remaining special cases are treated by explicit geometric considerations in Section \ref{special cases}.

Finally, in Section \ref{quadric moduli}, we generalize the techniques used in Section  \ref{lkq} to compute the Brauer groups of moduli of quadrics.



\section{Recollections on Brauer groups}\label{brauer}

In this section, we collect some well-known results  on Brauer groups,
without any claim to originality.
 Throughout, we shall
work with sheaves and cohomology   for the \'etale site.
Suppose $X$ is a scheme, and $n\geq 1$ be an integer.
Denote by $\PGL_{n,X}$ the sheaf of groups   obtained
by sheafifying the presheaf $U\mapsto \PGL_n(\Gamma(U,\O_X))$.
The short exact sequence of group-valued sheaves
$$
0\lra \GG_{m,X}\lra \GL_{n,X} \lra \PGL_{n,X} \lra 1
$$
yields a coboundary map 
$$
H^1(X,\PGL_{n,X})\lra H^2(X,\GG_{m,X}),
$$
whose image consists of $n$-torsion elements.
The union of these images for all integers $n\geq 1$ is called the \emph{Brauer group} $\Br(X)\subset H^2(X,\GG_m)$.
This subset is actually a subgroup, and it is contained in the torsion part.
The latter is usually called the \emph{cohomological Brauer group} and denoted by $\Br'(X)\subset H^2(X,\GG_m)$.
By Gabber's result (see \cite{de Jong 2006}), the inclusion $\Br(X)\subset \Br'(X)$ is an equality
if $X$  is quasicompact, separated, and admits an ample invertible sheaf.

If $X$ is proper over the field of complex numbers, then $\Br'(X)$ is   a group of the
form $(\QQ/\ZZ)^{\oplus n}\oplus T$, where $n=b_2-\rho$ is the difference of the second Betti number of 
the associated complex-analytic space $X^\an$ and the Picard number, and $T$ is the torsion part of $H^3(X^\an,\ZZ)$, see
for example \cite{Schroeer 2005}, Section 1.

One may regard the   set $H^1(X,\PGL_{n,X})$ as the set of isomorphisms classes of
$\PGL_{n,X}$-torsors. Since the  homomorphism of presheaves $\PGL_{n}(\Gamma(U,\O_X))\ra\Aut(\PP^{n-1}_U)$
becomes bijective upon sheafification,  the isomorphism classes of $\PGL_{n,X}$-torsors
are in bijective correspondence to the isomorphism classes of twisted forms
of $\PP^{n-1}_X$. 
Let us call a morphism $P\ra X$ a \emph{Brauer--Severi scheme} of relative dimension $n-1$ 
if there is some \'etale surjection $U\ra X$ with $P_U\simeq\PP^{n-1}_U$.
If $X=\Spec(k)$ is the spectrum of a field, one refers to  $P$ also as a \emph{Brauer--Severi variety}.
In the general case, one also says that $P\ra X$ is a \emph{family of Brauer--Severi varieties},
or a \emph{projective bundle}.

Given a Brauer--Severi scheme $f:P\ra X$ of relative dimension $n-1$,
the  coboundary of its isomorphism class is written as $[P]\in\Br(X)\subset H^2(X,\GG_m)$ and
called the \emph{Brauer class}.
This is an element whose order is a divisor of $n$. One also has the following geometric interpretation:
Clearly, the map $\O_X\ra f_*(\O_P)$ is bijective, 
such that we have an identification $\GG_{m,X}=f_*(\GG_{m,P})$. Moreover,
the degree map gives an identification $R^1f_*(\GG_m)=\ZZ_X$. With these identifications,
the Leray--Serre spectral sequence
$E_2^{pq}=H^p(X,R^qf_*(\GG_m))\Rightarrow H^{p+q}(P,\GG_m)$ yields an exact sequence
\begin{equation}
\label{exact sequence}
0\lra\Pic(X)\lra\Pic(P)\lra H^0(X,\ZZ_X)\lra H^2(X,\GG_m)\lra H^2(P,\GG_m).
\end{equation}

\begin{proposition}
\mylabel{image differential}
The image of the constant section $1_X$ under the differential   $H^0(X,\ZZ)\ra H^2(X,\GG_m)$
coincides with the Brauer class $[P]\in\Br(X)$. Moreover, 
the induced map $\Br(X)\ra \Br(P)$ is surjective.
\end{proposition}

\proof
The first statement is due to Giraud \cite{Giraud 1971}, Chapter V, Theorem 4.8.3.
According to Gabber \cite{Gabber 1981}, Theorem 2 on p.\ 193, the induced map on cohomological Brauer groups
$\Br'(X)\ra\Br'(P)$ is surjective. To see that the map on Brauer groups is surjective as well,
let $Q\ra P$ be a projective bundle. Choose some cohomology class $\alpha\in\Br'(X)$ with $f^*(\alpha)=[Q]$.
The structure morphism $f:P\ra X$ is flat and projective, the latter because the determinant of the dual of $\Omega_{P/X}^1$ 
is relatively ample.
According to a result of Edidin et al.\ \cite{Edidin; Hassett; Kresch; Vistoli 1999}, Theorem 3.6, the class $\alpha$ comes from some projective bundle over $X$.
\qed

\medskip
Let $P\ra X$ be a Brauer--Severi scheme of relative dimension $n-1\geq 0$. A closed subscheme
$L\subset P$ is called \emph{linear subscheme} of relative dimension $r-1\geq 0$, if there is an \'etale surjection $U\ra X$,
a locally free $\O_U$-module $\shE$  of rank $n$, a locally free quotient $\shE\ra \shF$ of rank $r$,
and an isomorphism $\varphi:P_U\ra\PP(\shE)$ inducing an isomorphism $L_U\ra\PP(\shF)$.
Here $\PP(\shF)\subset\PP(\shE)$ is the closed embedding obtained by forming the relative projective
spectrum for the surjection of graded rings  $\Sym^\bullet(\shE)\ra\Sym^\bullet(\shF)$.
Over fields, the following fact is due to Artin \cite{Artin 1982}, Proposition 3.6, with some different argument:

\begin{proposition}
\mylabel{linear subscheme}
Let $P\ra X$ be a Brauer--Severi scheme, and $L\subset P$
a linear subscheme. Then $[L]=[P]$ as elements of $\Br(X)$.
\end{proposition}

\proof
Let $f:P\ra X$ be the structure morphism, and $g=f|L$ its restriction to the linear subscheme.
Since the degree of an invertible sheaf on $\PP^n_R$ over some ring $R$ is the same as the degree of
its restriction to some linear subscheme, we have a commutative diagram
$$
\begin{CD}
R^1f_*\GG_{m,P} @>\deg>> \ZZ_X\\
@V\res VV                     @VV\id V\\
R^1g_*\GG_{m,L} @>\deg>> \ZZ_X,
\end{CD}
$$
of abelian sheaves on $X$. It follows that the vertical map on the left is bijective.
Using the Leray--Serre spectral sequences (\ref{exact sequence}) for both $f$ and $g$, we obtain a commutative diagram
$$
\begin{CD}
\Pic(P) @>>> H^0(X,\ZZ_X) @>>> H^2(X,\GG_m) @>>> H^2(P,\GG_m)\\
@VVV         @VV\id V              @VV\id V @VVV\\
\Pic(L) @>>> H^0(X,\ZZ_X) @>>> H^2(X,\GG_m) @>>> H^2(L,\GG_m)
\end{CD}
$$
and the result follows from Proposition \ref{image differential}.
\qed

\begin{proposition}
\mylabel{restriction function field}
Suppose  $X$ is a integral normal noetherian scheme, with   function field $F=\kappa(X)$.
If the strictly local rings $\O_{X,\bar{a}}^s$
are factorial for all closed geometric points $\bar{a}:\Spec(\Omega)\ra X$, 
then the restriction map $\Br(X)\ra\Br(F)$ is injective.
\end{proposition}

\proof
Let $P\ra X$ be a   Brauer--Severi scheme of relative dimension $n-1$ whose Brauer class $[P]$ maps to zero in $\Br(F)$.
Then the generic fiber $P_F$ is isomorphic to $\PP^{n-1}_F$, and in particular there
is an effective Cartier divisor $D_F\subset P_F$ of degree one. Consider its schematic closure $D\subset F$.
We claim that the ideal $\shI\subset\O_P$ is invertible. This question is local with
respect to the \'etale topology, thus we may assume that $X=\Spec(R)$ is strictly local and $P=\PP^{n-1}_X$.
By assumption, the ring $R$, which is normal and integral, is factorial.
By the Gau\ss\ Lemma,   polynomial rings over $R$ remain factorial, such that $P$ is locally factorial.
In turn, the Weil divisor $D\subset X$ is Cartier, such that $\shI\subset \O_X$ is invertible.
Let $\shL=\O_P(D)$ be the dual invertible sheaf. Since $X$ is connected, we have $\deg(\shL)=1_X$ as
global section of $\ZZ_X$. In light of the exact sequence (\ref{exact sequence}), the image of $1_X$ in $H^2(X,\GG_m)$
vanishes. According to Proposition \ref{image differential},  the Brauer class $[P]$ vanishes.
\qed

\medskip
We now recall a purity result of Grothendieck \cite{GB} in the following form.
For more general versions, see  Gabber \cite{Gabber 1998}.

\begin{proposition}
\mylabel{restriction bijective}
Let $k$ be a field and $X$   a smooth quasiprojective $k$-scheme. Let $Z\subset X$ be a closed
subset of codimension $\geq 2$, and $U=X\smallsetminus Z$ the corresponding open 
subscheme. Then the restriction map $\Br(X)\ra\Br(U)$ yields a bijection between    $l$-torsion parts,
for every integer $l\geq 1$ that is invertible in $k$.
\end{proposition}

\proof
Injectivity follows from Proposition \ref{restriction function field}. Surjectivity on $l$-torsion parts of $\Br'(X)\ra\Br'(U)$
is ensured by  \cite{GB}, Theorem 6.1. These cohomological Brauer groups are in fact
Brauer groups, because $X$ is quasiprojective \cite{de Jong 2006}.
\qed

\begin{proposition}
\mylabel{pullback bijective}
Let $R$ be a   ring, and $m\geq 0$ an integer. Then the pullback map 
$\Br(R)\ra\Br(\PP^m_R)$ is bijective. If $R$ is  normal, integral, noetherian, and 
the strictly local rings $\O_{\AA^m_R,\bar{a}}^s$
are factorial for all closed geometric points $\bar{a}:\Spec(\Omega)\ra \AA^m_R$, 
 then  $\Br(R)\ra\Br(\AA^m_R)$ induces a bijection between $l$-torsion parts, for all integers $l\geq 1$
that are invertible in $R$.
\end{proposition}

\proof
Both pullback maps are injective, because the projections $\PP^m_R\ra\Spec(R)$
and $\AA^m_R\ra\Spec(R)$ admit sections. By Proposition \ref{image differential},
the map $\Br(R)\ra\Br(\PP^m_R)$ is also surjective.

Now suppose that $R$ is normal, integral, noetherian, and all strictly
local rings of the polynomial ring $A=R[T_1,\ldots,T_m]$ at maximal ideals are factorial.
Let $\alpha\in\Br(A)$ be a cohomology class.
Denote by  $R\subset F$   the field of fractions, and $\hat{A}=R[[T_1,\ldots,T_m]]$ the completion.
Applying \cite{DeMeyer 1975}, Corollary 4 inductively, we see that the inclusion $R\subset \hat{A}$ induces a bijection
on Brauer groups. 
 Thus there is a projective bundle $P\ra\Spec(R)$ 
so that the Brauer class of  $P\otimes_R \hat{A}$ in $\Br(\hat{A})$ is the pullback of $\alpha$.

In light of Proposition \ref{restriction function field}, we are done if 
the Brauer class $[P]$ and the cohomology class $\alpha$ become equal in $\Br(A\otimes_RF)$. In particular, we may assume that $R=F$ is
a field. Now we can argue as follows: According to \cite{Auslander; Goldman 1960}, Proposition 7.6
the cokernel of the injection $\Br(K)\ra\Br(K[T])$  is $p$-torsion, where $K$ is any field and $p$ is its characteristic exponent.
Using induction on $m\geq1$, together  with Proposition \ref{restriction  function field}, we deduce that
the cokernel of the injection $\Br(F)\subset\Br(A)$ is $p$-torsion.
In turn, the $l$-torsion class $\alpha$ lies in the image of $\Br(F)$, and whence coincides with
the preimage of $[P]$.
\qed
\section{Recollections on quiver representations and their moduli}\label{quiver}

We recall all necessary facts on quiver representations, their moduli spaces and in particular smooth modes; for further details, the reader is referred to \cite{moduli}. For simplicity, we work over a ground field $k$ which is assumed to be algebraically closed of characteristic $0$.

Let $Q$ be a quiver with finite set of vertices $Q_0$ and finite set of arrows $Q_1$. An arrow $\alpha\in Q_1$ from a vertex $i\in Q_0$ to a vertex $j\in Q_0$ will be denoted by $\alpha:i\rightarrow j$. Denote by $\Lambda={\bf Z}Q_0$ the standard lattice in the real vector space $\mathbb{R}Q_0$, and by $\Lambda^+=\mathbb{N}Q_0$ the standard cone in $\Lambda$. The elements of $\Lambda^+$ will be called {\it dimension vectors} and are denoted by ${\bf d}=(d_i)_{i\in Q_0}$. Define the {\it Euler form} of $Q$, a bilinear form $\langle\_,\_\rangle$ on $\Lambda$, by
$$\langle {\bf d},{\bf e}\rangle=\sum_{i\in Q_0}d_ie_i-\sum_{\alpha:i\rightarrow j}d_ie_j.$$

A functional $\Theta\in(\mathbb{R}Q_0)^*$ is called a {\it stability}; it induces a {\it slope function} $\mu:\Lambda^+\smallsetminus\{0\}\rightarrow\mathbb{R}$ given by $\mu({\bf d})=\Theta({\bf d})/\dim {\bf d}$, where $\dim {\bf d}=\sum_id_i$.

A representation $V$ of $Q$ over an algebraically closed base field $k$ is given by (finite dimensional) $k$-vector spaces $V_i$ for $i\in Q_0$ and $k$-linear maps $V_\alpha:V_i\rightarrow V_j$ for $(\alpha:i\rightarrow j)\in Q_1$. The $k$-linear abelian category of all $k$-representations of $Q$ is denoted by ${\rm rep}_kQ$. It is of homological dimension at most one, that is, all ${\rm Ext}^i(\_,\_)$ for $i\geq 2$ vanish identically. There is a well-defined map ${\bf \dim}$ from the Grothendieck group $K_0({\rm rep}_kQ)$ to $\Lambda$ given by ${\rm\bf dim} V=(\dim V_i)_{i\in Q_0}$. Then, for all representations $V$ and $W$, we have $$\dim{\rm Hom}(V,W)-\dim{\rm Ext}^1(V,W)=\langle{\rm\bf dim}V,{\rm\bf dim} W\rangle.$$

The {\it slope} of a non-zero representation $V$ is defined as $\mu(V)=\mu({\rm\bf dim} V)$. The representation $V$ is called {\it semistable
(resp.~stable)} if $\mu(U)\leq\mu(V)$ (resp.~$\mu(U)<\mu(V)$) for all non-zero proper subrepresentations $U\subset V$. It is called {\it polystable} if it is isomorphic to a direct sum of stables of the same slope. The semistable representations of a fixed slope $\mu$ form an abelian subcategory ${\rm rep}^\mu_kQ\subset {\rm rep}_kQ$, whose simple (resp.~semisimple) objects are precisely the stables (resp.~polystables) of slope $\mu$. In particular, all stables have trivial endomorphism ring $k$. We have ${\rm Hom}({\rm rep}^\mu_kQ,{\rm rep}^\nu_kQ)=0$ if $\mu>\nu$.  Moreover, every representation $V$ admits a unique filtration (the
Harder--Narasimhan filtration)
$$0=V_0\subset V_1\subset\ldots\subset V_s=V$$
such that each subquotient $V_k/V_{k-1}$ is semistable, and
$$\mu(V_1/V_0)>\mu(V_2/V_1)>\ldots>\mu(V_s/V_{s-1}).$$

Given a dimension vector ${\bf d}=\sum_id_ii$ and $k$-vector spaces $V_i$ of dimension $d_i$ for $i\in Q_0$, we define the {\it variety of representations} $R_{\bf d}(Q)$ to be the affine scheme whose $k$-valued points are $\bigoplus_{\alpha:i\rightarrow j}{\rm Hom}_k(V_i,V_j)$. The $k$-points of $R_{\bf d}(Q)$ thus parametrize representations of $Q$ on the vector spaces $V_i$. The reductive algebraic group $G_{\bf d}=\prod_{i\in Q_0}{\rm GL}(V_i)$ acts on $R_{\bf d}(Q)$ via the base change action
$$(g_i)_i\cdot(V_\alpha)_\alpha=(g_jV_\alpha g_i^{-1})_{\alpha:i\rightarrow j}).$$

By definition, the orbits of $G_{\bf d}$ in $R_{\bf d}(Q)$ correspond to the isomorphism classes of $k$-representations of $Q$ of dimension vector ${\bf d}$. We denote by $R_{\bf d}^{\rm sst}(Q)$ (resp.~by $R_{\bf d}^{\rm st}(Q)$) the locus of semistable (resp.~stable) representations in $R_{\bf d}(Q)$. Then we have a chain of open and $G_{\bf d}$-stable embeddings
$$R_{\bf d}^{\rm st}(Q)\subset R_{\bf d}^{\rm sst}(Q)\subset R_{\bf d}(Q).$$
The quotient $M_{\bf d}^{\rm pst}(Q)$ of $R_{\bf d}^{\rm sst}(Q)$ by $G_{\bf d}$  exists. Its $k$-points parametrize isomorphism classes of polystable $k$-representations of $Q$ of dimension vector ${\bf d}$. Moreover, the geometric quotient of $R_{\bf d}^{\rm st}(Q)$ by $G_{\bf d}$ exists and is denoted by $M_{\bf d}^{\rm st}(Q)$; we denote by $q:R_{\bf d}^{\rm st}(Q)\rightarrow M_{\bf d}^{\rm st}(Q)$ the quotient map. The scheme $M_{\bf d}^{\rm st}(Q)$ admits an open embedding into $M_{\bf d}^{\rm pst}(Q)$, and its $k$-points parametrize the isomorphism classes of stable $k$-representations of $Q$ of dimension vector ${\bf d}$. If non-empty, the variety $M_{\bf d}^{\rm st}(Q)$ is smooth and irreducible of dimension $1-\langle{\bf d},{\bf d}\rangle$. In contrast, the variety $M_{\bf d}^{\rm pst}(Q)$ is typically singular (but still irreducible). Since the stabilizer of a point in $R_{\bf d}^{\rm st}(Q)$ reduces to the scalars, we have $M_{\bf d}^{\rm st}(Q)\simeq R_{\bf d}^{\rm st}(Q)/PG_{\bf d}$, where $PG_{\bf d}$ denotes the quotient of $G_{\bf d}$ by the subgroup $\Delta\simeq\mathbb{G}_m$, which is the image of the multiplicative group $\mathbb{G}_m$ embedded diagonally via $t\mapsto(t\cdot{\rm id}_{V_i})_i§$.\\[1ex]
We denote $g({\bf d}):=\gcd(d_i\, |\, i\in Q_0)$. If $g({\bf d})=1$, we can construct universal bundles $\mathcal{E}_i$ for $i\in Q_0$ on $M_{\bf d}^{\rm st}(Q)$ as follows: we consider the trivial bundle $V_i$ on $R_{\bf d}^{\rm st}(Q)$, and twist its natural $G_{\bf d}$-action by $(g_i)_i*v:=\chi((g_i)_i)^{-1}g_iv$, where $\chi\in X(G_{\bf d})$ denotes the character $\chi((g_i)_i)=\prod_{i\in Q_0}\det(g_i)^{a_i}$ for a choice of integers $a_i$ such that $\sum_ia_id_i=1$. Then the stabilizer of a point in $R_{\bf d}^{\rm st}(Q)$, which is the subgroup $\Delta$ defined above, acts trivially on $V_i$, and thus $V_i$ descends to a vector bundle $\mathcal{E}_i$ on $M_{\bf d}^{\rm st}(Q)$.

Given another dimension vector ${\bf n}=(n_i)_i\in\Lambda^+$ and $k$-vector spaces $W_i$ of dimension $n_i$ for $i\in Q_0$, there exists a scheme $M_{{\bf d},{\bf n}}^\Theta(Q)$ parametrizing pairs $(V,f)$ consisting of a semistable representation $V$ of $Q$ of dimension vector ${\bf d}$ and a tuple of maps $f=(f_i:W_i\rightarrow V_i)_{i\in Q_0}$ such that the following holds: if $U\subset V$ is a proper subrepresentation of $V$ such that $f_i(W_i)\subset U_i$ for all $i\in Q_0$, then $\mu(U)<\mu(V)$. Such pairs are parametrized up to the following equivalence relation: two such pairs $(V,f)$ and $(V',f')$ are equivalent if there exists an isomorphism $\varphi:V\rightarrow V'$ (given by linear maps $\varphi_i:V_i\rightarrow V_i'$ for all $i\in Q_0$) such that $f_i'=\varphi_if_i$ for all $i\in Q_0$.\\[1ex]
There exists a projective map $\pi_{\bf d}:M_{{\bf d},{\bf n}}^\Theta(Q)\rightarrow M_{\bf d}^{\rm pst}(Q)$, whose fiber over the stable locus is isomorphic to the projective space of dimension ${\bf n}\cdot{\bf d}-1=\sum_{i\in Q_0}n_id_i-1$. In fact, the restriction of $\pi_{\bf d}$ to $M_{\bf d}^{\rm st}(Q)$ defines a projective bundle $P_{\bf n}\rightarrow M_{\bf d}^{\rm st}(Q)$ of relative dimension ${\bf n}\cdot{\bf d}-1$, whose pullback to $R_{\bf d}^{\rm st}(Q)$ equals $\PP(\bigoplus_{i\in Q_0}V_i^{n_i})$. Therefore, if ${\bf n}'\leq{\bf n}$ componentwise, there exists a linear embedding $P_{{\bf n}'}\subset P_{\bf n}$.

In the special case where ${\bf n}$ has a single non-zero entry $n_i=1$, we also write $P_i=P_{\bf n}$.

\section{The main result}\label{main}

The aim of this section is to derive an explicit description of the Brauer group of a moduli space of stable quiver representations following closely the strategy of \cite{Balaji; Biswas; Gabber; Nagaraj 2007}.

Recall that a dimension vector ${\bf d}$ for a quiver $Q$ with stability $\Theta$ is called stable if $R_{\bf d}^{\rm st}(Q)\not=\emptyset$. For our purposes, the following notion will be useful.

\begin{definition} Let us a call a dimension vector ${\bf d}$ {\it amply stable} if 
$${\rm codim}_{R_{\bf d}(Q)}(R_{\bf d}(Q)\smallsetminus R_{\bf d}^{\rm st}(Q))\geq 2.$$
\end{definition}

\begin{theorem}\label{main theorem} Suppose ${\bf d}$ is an amply stable dimension vector.
Then the Brauer group ${\rm Br}(M_{\bf d}^{\rm st}(Q))$ is cyclic of order $g({\bf d})=\gcd(d_i\, |\, i\in Q_0)$, and the class of every $P_{\bf n}$ for ${\bf n}\not=0$ is a generator.
\end{theorem}

\proof Inside the projective bundle $P_{\bf d}\subset M_{{\bf d},{\bf d}}^\Theta(Q)$, we consider the open subset $U$ of equivalence classes of pairs $(V,f)$ such that every $f_i$ is an isomorphism and $V$ is stable. By definition of the smooth models, $U$ is isomorphic to $R_{\bf d}^{\rm st}(Q)$. We consider the following maps of Brauer groups:
$${\rm Br}(P_{\bf d})\longrightarrow{\rm Br}(U)\simeq{\rm Br}(R_{\bf d}^{\rm st}(Q))\longleftarrow{\rm Br}(R_{\bf d}(Q)).$$
The first map is injective by Proposition \ref{restriction function field}, and the third map is an isomorphism by ample stabilty of ${\bf d}$ and Proposition \ref{restriction bijective}. But $R_{\bf d}(Q)$ is just an affine space, thus its Brauer group is trivial by Proposition \ref{pullback bijective}, and we conclude that ${\rm Br}(P_{\bf d})$ is trivial, too. Now the exact sequence
$$\mathbb{Z}[P_{\bf d}]\longrightarrow{\rm Br}(M_{\bf d}^{\rm st}(Q))\longrightarrow{\rm Br}(P_{\bf d})$$
resulting from Proposition \ref{image differential} shows that ${\rm Br}(M_{\bf d}^{\rm st}(Q)$ is cyclic with generator $[P_{\bf d}]$.

Since there exist linear embeddings $P_{{\bf n}'}\rightarrow P_{\bf n}$ whenever ${\bf n}'\leq{\bf n}$ componentwise, we see that $\gamma=[P_{\bf n}]=[P_{\bf d}]$ for every ${\bf n}\not=0$ by Proposition \ref{linear subscheme}. Since $P_{\bf n}$ is of relative dimension ${\bf n}\cdot{\bf d}-1$, we infer from the discussion preceding Proposition \ref{image differential} that the order of $\gamma$ is necessarily a divisor ${\bf n}\cdot{\bf d}$ for all ${\bf n}\not=0$, and thus of $g({\bf d})$.

Suppose that the order of $\gamma$ equals some $0<h<g({\bf d})$. Then, for every $i\in Q_0$, the class of $\bigwedge^hP_{i}$ is trivial in ${\rm Br}(M_{\bf d}^{\rm st})$, thus $\bigwedge^hP_{i}$ is the projectivization of a vector bundle $W$ on $M_{\bf d}^{\rm st}(Q)$. The pullback $\hat{W}$ of $W$ to $R_{\bf d}^{\rm st}(Q)$ is therefore a $G_{\bf d}$-equivariant bundle (with the subgroup $\Delta\subset G_{\bf d}$ acting trivially), such that $\PP(\hat{W})\simeq \PP(\bigwedge^hV_i)$ $PG_{\bf d}$-equivariantly on $R_{\bf d}^{\rm st}(Q)$; here we regard $V_i$ as a trivial bundle on $R_{\bf d}(Q)$ and on the open subset $R_{\bf d}^{\rm st}(Q)$. Isomorphism of these projective bundles implies that there exists a $G_{\bf d}$-linearized line bundle $L$ on $R_{\bf d}^{\rm st}(Q)$ such that $\hat{W}\simeq\bigwedge^hV_i\otimes L$ $G_{\bf d}$-equivariantly, thus $\Delta$ acts trivially on $\bigwedge^h V_i\otimes L$. Since $R_{\bf d}^{\rm st}(Q)$ is open in the affine space $R_{\bf d}(Q)$, the line bundle $L$ is trivial, with
$G_{\bf d}$-linearization given by a character $(g_i)_i\mapsto \prod_i\det(g_i)^{c_i}$ of $G_{\bf d}$. The action of $\lambda\in\Delta$ on $\bigwedge^hP_i\otimes L$ is multiplication by the $(h+\sum_ic_id_i)$-th power of $\lambda$. By triviality of the action of $\Delta$ on $\bigwedge^hV_i\otimes L$, we conclude that $h=-\sum_ic_id_i$. Therefore $g({\bf d})$ divides $-\sum_ic_id_i=h$. This contradicts $0<h<g({\bf d})$.\qed

\medskip In the next section, we will derive a sufficient criterion for ample stability of a dimension vector, which is strong enough to allow computation of all Brauer groups in the case of multiple loop and generalized Kronecker quivers in Section \ref{lkq}. However, experiments suggest that the condition of ${\bf d}$ being amply stable is not essential. Therefore, we formulate:

\begin{conjecture}\label{main conj} If ${\bf d}$ is a stable dimension vector, the Brauer group ${\rm Br}(M_{\bf d}^{\rm sst}(Q))$ is cyclic of order $g({\bf d})$, and the class of every $P_{\bf n}$ for ${\bf n}\not=0$ is a generator.
\end{conjecture}

We consider the {\it tautological quiver representation} $\mathcal{V}$ on $R_{\bf d}^{\rm st}(Q)$, which is a representation of $Q$ in the category of locally free sheaves on $R_{\bf d}^{\rm st}(Q)$, such that $\mathcal{V}_i$ equals the constant sheaf corresponding to $V_i$ on $R_{\bf d}^{\rm st}(Q)$, see \cite{King 1994}. 

\begin{theorem} Let $Q$ be a quiver, $\Theta$ a stability and ${\bf d}$ a dimension vector such that the previous conjecture holds. Suppose that $g({\bf d})\geq 2$. Then there is no representation $\mathcal{E}$ of $Q$ into locally free sheaves on $M_{\bf d}^{\rm st}(Q)$ such that the pullback $q^*\mathcal{E}$ along the quotient map is $G_{\bf d}$-equivariantly isomorphic to $\mathcal{V}\otimes\mathcal{L}$ for some $G_{\bf d}$-linearized invertible sheaf $\mathcal{L}$ on $R_{\bf d}^{\rm st}(Q)$.
\end{theorem}

\proof Suppose, to the contrary, that such a quiver representation $\mathcal{E}$ exists. Then there exists a $G_{\bf d}$-equivariant isomorphism $$q^*\PP(\mathcal{E}_i)\longrightarrow \PP(\mathcal{V}_i)=q^*P_i.$$
This isomorphism descends to an isomorphism $\PP(\mathcal{E}_i)\rightarrow P_i$, contradicting non-vanishing of the Brauer class of $P_i$.\qed

\medskip In particular, in this case there is no {\it universal} quiver representation on $M_{\bf d}^{\rm st}(Q)$, in other words, $M_{\bf d}^{\rm st}(Q)$ is not a fine moduli space. Moreover, there is no {\it tautological} quiver representation $\mathcal{E}$ on $M_{\bf d}^{\rm st}(Q)$, in the sense that $q^*\mathcal{E}$ is isomorphic to the tautological quiver representation $\mathcal{V}$ on $R_{\bf d}^{\rm st}(Q)$.

In light of Proposition \ref{restriction function field}, the same applies to any non-empty open subset $U\subset M_{\bf d}^{\rm st}(Q)$.

\section{Reduction to combinatorics}\label{combinatorics}

In this section we use stratifications of varieties of representations to give a sufficient, purely combinatorial criterion for the codimension condition of Theorem \ref{main theorem}.

\begin{proposition}\label{suff} Let ${\bf d}$ be a stable dimension vector and suppose that, for each proper decomposition ${\bf d}={\bf e}+{\bf f}$  with $\mu({\bf e})\geq\mu({\bf f})$, we have $\langle{\bf e},{\bf f}\rangle\leq -2$. Then ${\bf d}$ is amply stable.
\end{proposition}

\proof We write $R_{\bf d}(Q)\smallsetminus R_{\bf d}^{\rm st}(Q)$ as the union of $R_{\bf d}(Q)\smallsetminus R_{\bf d}^{\rm sst}(Q)$ and $
R_{\bf d}^{\rm sst}(Q)\smallsetminus R_{\bf d}^{\rm st}(Q)$ and prove the desired codimension estimate separately for both subsets.

We first study $R_{\bf d}(Q)\smallsetminus R_{\bf d}^{\rm sst}(Q)$ and recall the Harder--Narasimhan stratification of $R_{\bf d}(Q)$: for every decomposition ${\bf d}={\bf d}^1+\ldots+{\bf d}^s$ into nonzero dimension vectors ${\bf d}^k$ such that $\mu({\bf d}^1)>\ldots>\mu({\bf d}^s)$ and such that $R_{{\bf d}^k}^{\rm sst}(Q)\not=\emptyset$ for all $k$, denote by $R_{\bf d}^{{\bf d}^*}(Q)$ the locus of all representations whose
Harder--Narasimhan filtration $V_*$ has subquotients $V_k/V_{k-1}$ of dimension vector ${\bf d}^k$ for $k=1,\ldots,s$. Then every $R_{\bf d}^{{\bf d}^*}(Q)$ is locally closed, we have $R_{\bf d}^{({\bf d})}(Q)=R_{\bf d}^{\rm sst}(Q)$, and $R_{\bf d}(Q)$ is the disjoint union of all the Harder--Narasimhan strata. Moreover, it is known by \cite[Proposition 3.4]{Reineke 2003} that
$${\rm codim}_{R_{\bf d}(Q)}(R_{\bf d}^{{\bf d}^*}(Q))=-\sum_{k<l}\langle {\bf d}^k,{\bf d}^l\rangle.$$

We will use this formula to prove that, under the hypothesis of the proposition, we have ${\rm codim}_{R_{\bf d}(Q)}(R_{\bf d}(Q)\smallsetminus R_{\bf d}^{\rm sst}(Q))\geq 2$. Namely, suppose that there exists a proper Harder--Narasimhan stratum of codimension $1$, thus a decomposition ${\bf d}^*$ such that $\sum_{k<l}\langle{\bf d}^k,{\bf d}^l\rangle=-1$. For every pair of indices $k<l$, we have $\langle{\bf d}^l,{\bf d}^l\rangle\leq 0$; namely, since $R_{{\bf d}^k}^{\rm sst}(Q),R_{{\bf d}^l}^{\rm sst}(Q)\not=\emptyset$, we can choose semistable representations $V$ and $W$ of dimension vector ${\bf d}^k$ and ${\bf d}^l$, respectively. Since $\mu(V)>\mu(W)$ by assumption, we have ${\rm Hom}(V,W)=0$. Thus, there exists precisely one pair $k_0<l_0$ such that $\langle{\bf d}^{k_0},{\bf d}^{l_0}\rangle=-1$, and $\langle {\bf d}^k,{\bf d}^l\rangle=0$ for all other pairs $k<l$. Defining ${\bf e}={\bf d}^1+\ldots+{\bf d}^{k_0}$ and ${\bf f}={\bf d}^{k_0+1}+\ldots+{\bf d}^s$, we thus have $\langle {\bf e},{\bf f}\rangle=-1$ and $\mu({\bf e})>\mu({\bf f})$. This contradicts our assumption.

Now we turn to $R_{\bf d}^{\rm sst}(Q)\smallsetminus R_{\bf d}^{\rm st}(Q)$.
Fix a decomposition ${\bf d}={\bf e}+{\bf f}$ into non-zero ${\bf e},{\bf f}$ such that $\mu({\bf e})=\mu({\bf d})=\mu({\bf f})$. We define ${\rm Grass}_{\bf e}({\bf d})$ as the product of Grassmannians $\prod_{i\in Q_0}{\rm Grass}_{e_i}(V_i)$. Inside the scheme $R_{\bf d}^{\rm sst}(Q)\times{\rm Grass}_{\bf e}({\bf d})$, we consider the closed subscheme $X_{{\bf e},{\bf f}}$ of pairs $((V_\alpha)_\alpha,(U_i)_i)$ such that $V_\alpha(U_i)\subset U_j$ for all arrows $\alpha:i\rightarrow j$ of $Q$. In other words, $X_{{\bf e},{\bf f}}$ parametrizes semistable representations together with a subrepresentation of dimension vector ${\bf e}$. The projection $p_1:X_{{\bf e},{\bf f}}\rightarrow R_{\bf d}^{\rm sst}(Q)$ is projective, with image consisting of all semistable representations admitting a subrepresentation of dimension vector ${\bf e}$. By definition, we thus have
$$R_{\bf d}^{\rm sst}(Q)\smallsetminus R_{\bf d}^{\rm st}(Q)=\bigcup_{{\bf e},{\bf f}}p_1(X_{{\bf e},{\bf f}}),$$
where the sum ranges over all proper decompositions ${\bf d}={\bf e}+{\bf f}$ into dimension vectors of the same slope.
The projection $p_2:X_{{\bf e},{\bf f}}\rightarrow{\rm Grass}_{\bf e}({\bf d})$ is a homogeneous $G_{\bf d}$-bundle, whose relative dimension $r$ is easily computed (see e.g. \cite{CFR1}) as
$$r=\sum_{\alpha:i\rightarrow j}(e_ie_j+f_ie_j+f_if_j).$$
A standard Euler form calculation then allows us to estimate
$${\rm codim}_{R_{\bf d}^{\rm sst}(Q)}(p_1(X_{{\bf e},{\bf f}}))\geq \dim R_{\bf d}(Q)-\dim {\rm Grass}_{\bf e}({\bf d})-r=-\langle{\bf e},{\bf f}\rangle\geq 2$$
This gives the desired codimension estimate.\qed



\section{Loop and Kronecker quivers}\label{lkq}

We now show that Theorem \ref{main theorem} and the combinatorial criterion Proposition \ref{suff}, combined with two explicit calculations, suffice to compute the Brauer groups of all non-trivial moduli spaces for multiple loop and generalized Kronecker quivers.

Let $L_m$ be the {\it $m$-loop quiver} with a single vertex and $m\geq 0$ loops; let $K_m$  be the {\it $m$-arrow Kronecker quiver} with two vertices $1$ and $2$ and $m$ arrows from $1$ to $2$.

\begin{theorem}\label{lk} Conjecture \ref{main conj} holds for multiple loop quivers $L_m$ and generalized Kronecker quivers $K_m$.
\end{theorem}

The proof will occupy this and the following section.

We start with the quiver $L_m$. For $m=0,1$, the only non-trivial moduli spaces of stable representations occur in dimension one, and are isomorphic to a point, resp.~an affine line. By abuse of notation, we identify a dimension vector ${\bf d}=(d)$ with the integer $d$. The moduli space $M_{d}^{\rm st}(L_m)$ parametrizes $m$-tuples of $d\times d$-matrices without nontrivial common invariant subspaces up to simultaneous conjugation. We shall use Theorem \ref{main theorem} and Proposition \ref{suff}. Assume that the criterion of Proposition \ref{suff} is not fulfilled. Then $d=e+f$ with $e,f\geq 1$ and  $(m-1)ef\leq1$, which holds if and only if $m=2$, $e=f=1$, and thus $d=2$.

In this case, the relevant moduli space $X=M_2^{\rm st}(L_2)$ is well known to be isomorphic to the open subset of $\mathbb{A}^5$ (with coordinates $a,b,c,d,e$) given by $b^2\not=ac$. Namely, the coordinates $a,\ldots,e$ correspond, respectively,  to the invariants
$${\rm tr}(A'^2), {\rm tr}(A'B'), {\rm tr}(B'^2), {\rm tr}(A), {\rm tr}(B)$$
of two $2\times 2$-matrices $A$ and $B$, where $A'=A-\frac{1}{2}{\rm tr}(A)E$, $B'=B-\frac{1}{2}{\rm tr}(B)E$. As will be proved in the next section, this description allows to identify the Brauer group of $X$ as $\ZZ/2\ZZ$. Moreover, the bundle $P_1\rightarrow X$ is given by the equation
$$cx^2+az^2=2(y^2+bxz)$$
in $X\times\mathbb{P}^2$, with homogeneous coordinates $(x:y:z)$ for $\mathbb{P}^2$. Namely, the coordinates $x,y,z$ correspond, respectively, to the semiinvariants
$$\det(v|Av), \det(v|Bv), \det(Av|Bv)$$
for the natural ${\rm GL}_2(k)$-action on triples $(A,B,v)$ of two $2\times 2$-matrices and a vector in $k^2$. As will be proved in the next section, the class of this bundle in the Brauer group is non-vanishing. We have thus proved that ${\rm Br}(M_d^{\rm st}(L_m))\simeq\mathbb{Z}/d\mathbb{Z}$, and the class of any $P_n$ for $n\geq 1$ is a generator.

\bigskip Next we consider the generalized Kronecker quiver $K_m$.  We choose the stability $\Theta(d_1,d_2)=d_1$, which is the only relevant one, see \cite{moduli}. We then have:

\begin{proposition} Let $(d_1,d_2)$ be a stable dimension vector for $K_m$. Then $(d_1,d_2)$ is amply stable, except for the case $m=3$, $d_1=d_2=2$.
\end{proposition}

\proof First, we can assume that $m\geq 3$. Namely, in case $m=0$, the stable dimension vectors are $(1,0)$, $(0,1)$, in case $m=1$, they are $(1,0)$, $(0,1)$ and $(1,1)$, in case $m=2$ they are $(n,n+1)$,$(n+1,n)$, $(1,1)$ for $n\geq 0$. In all these cases, the moduli spaces are single points or (in the last case) a projective line.

Next, using reflection functors and duality (see \cite[Proposition 4.3]{Weist}), we can assume without loss of generality that
\begin{equation}\label{f0}d_1\leq d_2\leq \frac{m}{2}d_1.\end{equation} We write $(d_1,d_2)=(np,nq)$ for $n\geq 1$ and coprime $p$ and $q$. We can assume $mpq-p^2-q^2\geq 0$ (otherwise the moduli space of stables is empty or reduces to a single point).

Let us assume that we have a decomposition $(d_1,d_2)=(a,b)+(c,d)$ contradicting the hypothesis of Proposition \ref{suff}. Obviously $b=nq-d$ and $c=np-a$, and in particular
\begin{equation}\label{f1} a\leq np\mbox{ and }d\leq nq.\end{equation}
In light of the stability $\Theta(d_1,d_2)=d_1$, the slope condition $\mu(a,b)\geq\mu(c,d)$ is equivalent to $a/b\geq c/d$, which in turn means
\begin{equation}\label{f2}k:=pd+qa-npq\geq 0.\end{equation}
The Euler form condition $\langle(a,b),(c,d)\rangle\geq-1$, after some rewriting, reads
\begin{equation}\label{f3}pq\geq(mpq-p^2-q^2)ad+(pa+qd)k.\end{equation}
Our assumption (\ref{f0}) on the dimension vector implies
\begin{equation}\label{f4}p\leq q\leq\frac{m}{2}p.\end{equation}
If $a=0$, then $c=0$ by the slope condition, thus $p=0$, thus $q=0$, yielding the dimension vector ${\bf d}=0$ a contradiction, and similarly for the assumption $d=0$. Thus we can assume
$ad\geq 1$, and (\ref{f3}) yields the estimate $pq\geq mpq-p^2-q^2$,
which, after dividing by $pq$, reads
$$m\leq\frac{p}{q}+\frac{q}{p}+1.$$
Let us first treat the case $p=1=q$. Inserting in (\ref{f3}), we get $1\geq(m-2)ad+(a+d)k$, which by $a,d\geq 1$ implies $m=3$, $k=0$, $a=1=d$, and thus ${\bf d}=(2,2)$, as claimed in the proposition.

So let us now assume $(p,q)\not=(1,1)$, thus $p<q$ by coprimality.
We use (\ref{f4}) to estimate $\frac{p}{q}<1$ and $\frac{q}{p}\leq\frac{m}{2}$, thus
$$m<\frac{m}{2}+2,$$
which implies $m<4$ and thus $m=3$.
 Then (\ref{f4}) yields the estimate 
$$3pq-p^2-q^2=(3p-q)q-p^2\geq\frac{3}{2}pq-p^2>\frac{1}{2}p^2,$$
and thus by (\ref{f4}) and (\ref{f3}):
$$\frac{3}{2}p^2\geq pq\geq(3pq-p^2-q^2)ad>\frac{1}{2}p^2ad.$$
This yields $ad<3$, and thus $ad\in\{1,2\}$. This leaves us with the three cases $a=d=1$ or $a=1$, $d=2$ or $a=2$, $d=1$.

To estimate $n$, we use equation (\ref{f2}), which gives us
$$n=\frac{1}{q}d+\frac{1}{p}a-\frac{k}{pq}\leq\frac{1}{q}d+\frac{1}{p}a\leq d+a\leq 3,$$
and thus $n\leq 3$. In case $n=3$, the previous inequality yields $p=q=1$, a contradiction.

Thus we are in the special situation $m=3$ and $n,ad\in\{1,2\}$. Rewriting again the Euler form condition, we have
\begin{equation}\label{finfty}nap+ndq=a^2+d^2+3ad-1.\end{equation}
In each of the six cases $n=1,2$ and $(a,d)=(1,1),(1,2),(2,1)$, one sees directly that no pair $(p,q)$ satisfies (\ref{f2}), (\ref{f4}) and (\ref{finfty}) simultaneously.\qed


\medskip In light of Proposition \ref{suff}, the preceding result finishes the proof of Theorem \ref{lk} for generalized Kronecker quivers $K_m$, except for the remaining special case ${\bf d}=(2,2)$ for the quiver $K_3$. In this case, the moduli space $X=M_{(2,2)}^{\rm st}(K_3)$ is isomorphic to the open subset $X\subset \mathbb{P}^5$, with homogeneous coordinates $(a:b:c:d:e:f)$, such that
$$4adf+bce-c^2d-ae^2-b^2f\not=0.$$
These six coordinates essentially result from polarization of the determinant, that is, expressing $\det(\alpha A+\beta B+\gamma C)$ for a triple $(A,B,C)$ of $2\times 2$-matrices as a quadratic form in $\alpha,\beta,\gamma$. The bundle $P_1$ is the $\mathbb{P}^1$-bundle over this space given by the equation
$$fx^2-exy+dxz+cy^2-byz+az^2=0$$
in $X\times\mathbb{P}^2$ (with coordinates $(x:y:z)$ for $\mathbb{P}^2$). Similar to the case above, the coordinates $x,y,z$ correspond, respectively, to the determinants
$$\det(Av|Bv),\det(Av|Cv),\det(Bv|Cv)$$
for a quadruple $(A,B,C,v)$ with $A,B,C$ as above and $v\in k^2$. We will prove in the next section that the Brauer group of $X$ has two elements, and that the class of $P_i$ is non-vanishing.


\section{The two special cases}
\label{special cases}

The task now is to compute the Brauer group in the two special cases $M_2^{\rm st}(L_2)$ and $M_{(2,2)}^{\rm st}(K_3)$.

We start with the Kronecker quiver $K_3$ with dimension vector $\mathbf{d}=(2,2)$.
This moduli space can be identified with the homogeneous space $\GL_3/\GO_3$, by regarding the latter as an open subset of $\PP^5$ as follows:
Choose six indeterminates $a,b,\ldots,f$  and let
$$
h=\det
\begin{pmatrix} 
2a & b  & c\\
b  & 2d & e\\
c  & b  & 2f
\end{pmatrix}
=2(4adf+bce-c^2d-ae^2-b^2f)
$$
be the determinant of the generic symmetric $3\times 3$-matrix, which is a homogeneous
polynomial of degree three. Then the homogeneous space $\GL_3/\GO_3$
becomes isomorphic to the open set $U=D_+(h)\subset\PP^5$.

\begin{proposition}
The group $\Br(U)$ is cyclic of order two.
\end{proposition}

\proof
Consider the complementary closed subscheme $Y=V_+(h)\subset\PP^5$, which is a cubic fourfold.
Let $Z=\Sing(Y)$ be its singular subscheme, which is defined by
the jacobian ideal $J=(h,\partial h/\partial a,\ldots,\partial h/\partial f)$.
A computation with Magma \cite{Magma} reveals that $Z$ is integral and 2-dimensional.
Hence $Y$ is normal, by Serre's Criterion.

Now fix an integer $n\geq 1$, and consider the regular locus  $Y_0=Y\smallsetminus Z$.
Combining \cite{Ford 1992},  Lemma 0.1 and Theorem 1.1, we obtain an identification of $n$-torsion groups
$$
\Br(U)_n= H^3_{Y_0}(\PP^5\smallsetminus Z,\mu_n) = H^1(Y_0,\ZZ/n\ZZ).
$$
The latter is the set of isomorphism classes of $\ZZ/n\ZZ$-torsors over the smooth scheme $Y_0$.
This group becomes isomorphic to $\Pic(Y_0)_n$, after choosing an isomorphism $\mu_n\simeq\ZZ/n\ZZ$,
in light of \cite{Raynaud 1970}, Proposition 6.2.1.
Since $Y$ is normal, we have an identification $\Pic(Y_0)=\Cl(Y)$ with the class group
of Weil divisors modulo linear equivalence.
Our task now is  to show that the class group $\Cl(Y)$ is cyclic of order two.

To achieve this, we pass from cubic fourfolds to cubic surfaces.
Let $H,H'\subset \PP^5$ be two   hyperplanes so that the iterated hyperplane section
$S=Y\cap H\cap H'$ becomes 2-dimensional, such that $S\subset H\cap H'\simeq\PP^3$ is a cubic surface. 
According to the Lefschetz-Theorems of Ravindra and Srinivas for class groups (\cite{Ravindra; Srinivas 2006}, Theorem 1), 
the restriction map $\Cl(Y)\ra\Cl(S)$ is bijective,
provided that the hyperplanes are general.
Consider first the special hyperplanes $H=V_+(b-c)$ and $H'=V_+(b-e)$.
Then the cubic surface $S\subset\PP^3$ is given by the homogeneous polynomial $h'=4adf+b^3-c^2d-ac^2-b^2f$.
A computation with Magma reveals that the singular locus $Z'=\Sing(S)$
consists of four   closed points, which are contained in $Z$.
It follows that the singular subscheme $Z=\Sing(Y)$ has degree $\deg(Z)\geq 4$ as closed subscheme 
$Z\subset\PP^5$. 

Now let $H,H'\subset\PP^5$ be two general hyperplanes. Then the finite scheme $Z\cap H\cap H'$ consists of
at least four points, and  Bertini tells us that  the cubic surface $S$ has singular locus $\Sing(S)=Z\cap H\cap H'$.
Now recall that normal cubic surfaces are classified: The result goes back to Schl\"afli \cite{Schlaefli 1864}, 
where treated in a modern way by Bruce and Wall \cite{Bruce; Wall 1979}, and further refined by Sakamaki \cite{Sakamaki 2010}.
It follows from this classification that the normal cubic surface with at least four singularities is unique
up to coordinate change, that there are precisely four  singularities, each a  rational double points of type $A_1$,
and that this cubic surface may by described by the homogeneous equation
$$
x_3(x_0x_2-x_1^2) - (x_0-x_1)(x_1-x_2)=0
$$
in the variables $x_0,x_1, x_2$, compare \cite{Sakamaki 2010}, Theorem 2.

Finally, consider the minimal resolution of singularities $\tilde{S}\ra S$. Then $\tilde{S}$ is
a {\it weak del Pezzo surface} of degree $K_{\tilde{S}}^2=K_S^2=3$. By the structure theory
of weak del Pezzo surfaces, $\tilde{S}$ is obtained from $\PP^2$
by blowing-up six points, some of which may be  infinitesimal near.
Using \cite{Sakamaki 2010}, Section 2, in particular paragraph 2.2.7,  we can make this completely explicit:
Let $C\subset \PP^2=\Proj(\CC[x_0,x_1,x_2])$ be the quadric curve given by $x_0x_2-x_1^2=0$.
Then $\tilde{S}\ra \PP^2$ is obtained by first blowing-up the the three points
$$
(0:0:1),(1:0:0),(1:1:1)\in\PP^2,
$$
which lie on the quadric curve $C$, followed by
a blowing-up of the three intersection points of the resulting exceptional divisors $E_1,E_2,E_3$ with the strict transform of $C$.
The four $(-2)$-curves $\tilde{C}_i\subset\tilde{S}$, $1\leq i\leq 4$ arise as the strict transforms of $C$ and $E_1,E_2,E_3$.
A simple computation reveals that their sum $\sum\tilde{C}_i$ is not primitive in $\Pic(\tilde{S})$, 
in fact, it is exactly divisible by 2. 
More precisely, the quotient $\Pic(\tilde{S})/H$
by the subgroup $H\subset\Pic(\tilde{S})$ generated by the $(-2)$-curves $\tilde{C}_i$ has torsion part of order two.
This can be also seen by considering the discriminants for the bilinear forms on  $H\subset H^{\perp\perp}$.
Using the identification $\Cl(S)=\Pic(\tilde{S})/H$, we conclude that the class group $\Cl(S)$ has order two.
\qed

\medskip
Now let $x,y,z$ be three further indeterminates, and
consider the relative quadric $P\subset\PP^2\times \PP^5$   defined by the bi-homogeneous equation
$$
fx^2-exy + dxz +cy^2-byz + az^2=0.
$$
We are mainly interested in the restriction $P_U\ra U$ to the open subset 
$U=\PP^5\smallsetminus Y=\GL_3/\GO_3$ considered above. As explained in the previous section, this actually is
the smooth model for the quiver moduli space.

\begin{proposition}
The projection $P_U\ra U$ is a Brauer--Severi scheme whose Brauer class generates $\Br(U)$.
\end{proposition}

\proof
The fibers of the projection are smooth, because the symmetric matrix corresponding
to the equation is invertible over $U\subset\PP^5$. In turn, $P_U\ra U$ is a Brauer--Severi scheme.
Since $U$ is smooth, and $\Br(U)$ has order two, it suffices to check that the generic fiber $P_\eta$
contains no rational point. This easily follows, because in the function field $F=\kappa(\eta)$,
the five elements $a/f,\ldots,e/f\in K$ form a transcendence base.
\qed

\medskip
Finally, we consider the quiver moduli space for the loop quiver $L_2$ with
dimension vector $\mathbf{d}=(2)$. This quiver moduli space may be regarded as 
an open open subscheme of $\AA^5$ as follows: Choose five indeterminates $a,\ldots,e$ 
and consider the polynomial $h=b^2-ac$.
Then the quiver moduli space may be regarded as the affine open subscheme $U=D(h)\subset\AA^5$.
Let $x,y,z$ be three further indeterminates, and consider the closed subscheme $P\subset\PP^2\times\AA^5$
defined by the homogeneous equation
$$
cx^2+az^2=2(y^2+bxz).
$$
The restriction $P_U\ra U$ is actually the smooth model for the quiver moduli.

\begin{proposition}
The Brauer group $\Br(U)$ is cyclic of order two, and the projection $P_U\ra U$
is a Brauer--Severi variety whose Brauer class generates $\Br(U)$.
\end{proposition}

\proof
According to \cite{Ford 1989}, Example 2 for Theorem 1, the Brauer group of the localization $R=\CC[a,b,c][1/(b^2-ac)]$
has order two. In light of Proposition \ref{pullback bijective}, the same holds for the affine scheme $U$.
The second assertion follows as in the preceding proof.
\qed

\section{Brauer groups for moduli of quadrics}
\label{quadric moduli}

In this section, we compute the Brauer group for the Hilbert moduli space
of smooth odd-dimensional quadrics. In the special case of 1-dimensional quadrics, this moduli
space coincides with the quiver moduli space for the generalized Kronecker quiver $K_3$. It turns out that the computation works over arbitrary ground rings.
We start by recalling some facts on quadrics.

Let $S$ be a scheme, and $n\geq 0$ be an integer. 
A closed subscheme $X\subset \PP^{n+1}_S$ is called a \emph{relative quadric}
if the structure morphism $X\ra S$ is flat, and for each point $s\in S$ the fiber
$X_s\subset\PP^{n+1}_s$ is an effective Cartier divisor of degree two.
The corresponding invertible sheaf is of the form $\O_{\PP^{n+1}_S}(X)\simeq \O_{\PP^{n+1}_S}(2)\otimes\pr^*(\shN)$ for some
invertible $\O_S$-module $\shN$, in light  of the decomposition
$$
\Pic(\PP^{n+1}_S) = \ZZ\O(1) \oplus \Pic(S).
$$
Suppose for the moment that $S=\Spec(R)$ is affine and that $\shN$ is trivial.
Then $X=V_+(b)$ for some section 
$$
b\in H^0(\PP^{n+1}_S,\O_{\PP^{n+1}_S}(2))=\Sym^2(E_R),
$$
where $E=\ZZ^{\oplus(n+2)}$.  
With respect to the standard basis $T_0,\ldots,T_{n+1}\in E$, we may regard $b=(b_{ij})_{0\leq i,j\leq n+1}$
as a symmetric matrix with entries in $R$. This correspond to the quadratic form
$$
Q:E^\vee_R\lra R,\quad Q(\sum\lambda_ie_i)=\sum b_{ij}\lambda_i\lambda_j,
$$
where $e_i\in E^\vee$ is the basis dual to the standard basis $T_i\in E$, and
the sum runs over all subsets $\left\{i,j\right\}\subset \left\{0,\ldots,n+1\right\}$
of cardinality one or  two, as explained in \cite{A 9}, \S3, No.\ 4, Proposition 2. 
Note that the associated
bilinear form satisfies $\Phi(e_i,e_j)=b_{ij}$ for $i\neq j$, and $\Phi(e_i,e_i)=2b_{ii}$.
Moreover,  $b$ and whence $Q$
are unique up to unique unit  $u\in R$.

Now let $S$ be again arbitrary.
As outlined in \cite{SGA 7b}, Expose XII, we can construct a sheaf of associative 
$\O_S$-algebras $\shA$
that locally comes from the even parts $A=\Cl_+(Q_R)$ of the \emph{Clifford algebras} attached
to the quadratic forms $Q_R:E^\vee_R\ra R$.
To this end, let $U_\alpha\subset S$ be the collection of all affine open subsets $U_\alpha=\Spec(R_\alpha)$
so that the invertible sheaves
$\O_{\PP^{n+1}_S}(X)$ and $\O_{\PP^{n+1}_S}(2)$ become isomorphic over the preimage of $U_\alpha$.
As above, choose $b_\alpha\in \Sym^2(E_{R_\alpha})$, which gives the quadratic form $Q_\alpha$,
and let $A_\alpha=\Cl_+(Q_\alpha)$ the resulting even part of the Clifford algebra $\Cl(Q_\alpha)$.
We refer to \cite{A 9}, \S9 for details on Clifford algebras.
On the overlaps $U_{\alpha\beta}=U_\alpha\cap U_\beta$,  one has $Q_\beta=u_{\alpha\beta}Q_\alpha$
for some unique section $u_{\alpha\beta}\in\Gamma(U_{\alpha\beta},\O_S^\times)$.
By this uniqueness, the resulting cochain $(u_{\alpha\beta})$ satisfies the cocycle condition.
Since the \emph{even part} of Clifford algebras is functorial with respect to similitudes, it follows that
there are unique isomorphisms
$$
[u_{\alpha\beta}]: \widetilde{\Cl}_+(Q_\alpha)|U_{\alpha\beta}\lra \widetilde{\Cl}_+(Q_\beta)|U_{\alpha\beta},
$$
of quasicoherent sheaves (confer \cite{Wonnenburger 1962} and \cite{SGA 7b}, Expose XII, Lemma 1.3.1), which also satisfy the cocycle condition.
In turn, we have a descend datum for the  the sheaves $\shA_\alpha=\widetilde{\Cl}_+(Q_\alpha)$, which
yields the desired sheaf of associative $\O_S$-algebras $\shA$. By construction,
the $\O_X$-module $\shA$ is locally free of rank $2^{n+1}$. Up to unique isomorphism,
it does not depend on the choice
of the local sections $b_\alpha$.

\begin{proposition}
\mylabel{smooth Azumaya}
If $n$ is odd, and the structure morphism $X\ra S$ of the relative quadric $X\subset\PP^{n+1}_S$ is smooth, 
then the $\O_S$-algebra $\shA$ is an Azymaya algebra.
\end{proposition}

\proof
According to \cite{GB}, Theorem 5.1, it suffices to treat the case that $X=\Spec(k)$ is the
spectrum of an algebraically closed field. Set $V=E_k^\vee$, and let $Q:V\ra k$ be a quadratic form defining
the quadric $X\subset\PP^{n+1}_k$.
Write $n+1=2r$. Suppose first that this quadratic form $Q$ equals  the  \emph{standard quadratic form} 
\begin{equation}
\label{standard quadric}
Q_\text{std}(\sum_{i=0}^{n+1}\lambda_iT_i)=\sum_{i=0}^{r-1}\lambda_i\lambda_{i+r} + \lambda_{n+1}^2.
\end{equation} 
For each $0\leq i\leq r-1$, let $V_i\subset V$  be the linear subspace generated by $e_i,e_{i+r}\in V$. 
Furthermore, let $V''\subset V$ be the linear subspace
generated bei $e_{n+1}\in V$. Write $V'=V_0\oplus\ldots\oplus V_{r-1}$, such that $V=V'\oplus V''$. Now consider
the restriction   $q',q_i$ of the quadratic form $Q$ to the subspaces
$V',V_i$, respectively. Then we have a decomposition of algebras
$$
\Cl_+(Q) = \Cl(q') = \Cl(q_0)\otimes\ldots\otimes\Cl(q_{r-1}),
$$
by \cite{Lam 2005}, Chapter V, Corollary 2.10.. 
Since the quadratic forms $q_i$ are nondegenerate, the Clifford algebras $\Cl(q_i)$ are Azumaya algebras,
by \cite{A 9}, \S9, No.\ 4, Theorem 2. In turn, the same holds for $\Cl_+(Q)$.

It remains to check that $Q$ is similiar to the standard quadratic form. Let $p\geq 0$ be the
characteristic of the field $k$. Consider first the case $p\neq 2$.
Since $X$ is smooth, the associcated symmetric bilinear form $\Phi$ is nondegenerate.
According to \cite{A 9}, \S4, No.\ 3, the quadratic form $Q$ is indeed similar to the standard form.
Finally, suppose that $p=2$. Then $\Phi$ is alternating, thus necessarily degenerate, because $\dim(V)=n+2$ is odd.
Let $V^\perp\subset V$ be the orthogonal complement of the whole vector space. Since $X$ is smooth,  
we may apply \cite{Buchweitz; Eisenbud; Herzog 1987}, Theorem 1.1
and deduce that $V^\perp$ is one-dimensonal and not \emph{singular} (in the sense of \cite{A 9}, \S4, No.\ 2).
Choose a linear complement $V=V^\perp\oplus V'$. Obviously, this is an orthogonal
complement, of even dimension, and  the restriction $Q|V'$ is nondegenerate. Applying \cite{A 9}, \S4, No.\ 3 to  the quadratic form
$Q|V'$, we infer that $Q$ is similar to the standard form.
\qed

\medskip
We say that the relative quadric $X\subset \PP^{n+1}_S$ is \emph{smooth} if
the structure morphism $X\ra S$ is smooth. 
Suppose this is the case and that  $n+2=2r+1\geq 3$ is odd.
The ensuing Azumaya algebra $\shA$ corresponds to  a Brauer--Severi scheme $B\ra S$
of relative dimension $2^r-1$. The corresponding Brauer class $[B]\in\Br(S)$ 
is called the \emph{Clifford invariant} of the smooth quadric $X\subset\PP^{n+1}_S$.

Let us now turn to the universal situation. Denote by $U'\subset\Hilb_{\PP^{n+1}}$ be
the open subscheme of the Hilbert scheme over the base scheme $S=\Spec(\ZZ)$
parameterizing closed subschemes $X\subset\PP^{n+1}_k$
with Hilbert polynomial
$$
\chi\O_X(t)=\chi\O_{\PP_k^{n+1}}(t) - \chi\O_{\PP_k^{n+1}}(t-2)= \binom{t+n+1}{n+1} - \binom{t+n-1}{n+1},
$$
and whose ideal sheaf is invertible. Clearly, $U'$ can be identified with projective
space attached to the graded ring $\Sym^\bullet(\Sym^2(E^\vee))$.
Then the universal closed subscheme restricted to $U'$ 
is a relative quadric.  Let $U\subset U'$ be the subset where the relative quadric
is smooth, which is open by \cite{EGA IVb}, Corollary 6.8.7. The restriction 
$X\subset \PP^{n+1}_U$ of the universal closed subscheme is called the \emph{universal smooth quadric}.
Note that the structure morphism $U\ra\Spec(\ZZ)$ is smooth and with geometrically integral fibers.
All fibers are nonempty, because
the standard quadratic form (\ref{standard quadric}) defines a section for the structure morphism.

The group scheme $\GL_{n+2}$ acts in the canonical way from the right on  $\PP^{n+1}$, and
induces an action  on the Hilbert scheme from the right. 
We thus obtain a morphism
$$
\GL_{n+2}\lra \Hilb_{\PP^{n+1}}, \quad A\longmapsto X_{\text{std}}\cdot A,
$$
where $X_{\text{std}}\subset\PP^{n+1}$ is the relative quadric defined by the standard quadratic form (\ref{standard quadric}).
Its stabilizer is the group scheme of \emph{orthogonal similitudes} $\GO_{n+2}$, and we thus get a morphism
$\GO_{n+2}\backslash\GL_{n+2}\ra \Hilb_{\PP^{n+1}}$, which factors
over the scheme of smooth quadrics $U\subset\Hilb_{\PP^{n+1}}$.
Note that the quotient actually exists as a scheme,  according to \cite{Anantharaman 1973}.

\begin{proposition}
\mylabel{identification moduli}
The induced morphism $\GO_{n+2}\backslash\GL_{n+2}\ra U$ is an isomorphism.
\end{proposition}

\proof
Set $V=\GO_{n+2}\backslash\GL_{n+2}$ and write $f:V\ra U$ for the morphism in question.
First, we verify that $f$ a universal homeomorphism.
Since $f$ is of finite presentation, it is enough to show that it is universally bijective, according to \cite{EGA IVa}, Corollary 1.10.4.
For this, it suffices to check
that for each algebraically closed field $k$, the induced morphism 
$f_k:V_k\otimes k\ra U_k$ is bijective. The latter is injective, because
the group of orthogonal similitudes is the stabilizer of the standard quadric.
It is surjective as well, because any smooth quadric comes from
a quadratic form that is similiar to the standard quadratic form, as we saw
in the proof for Proposition \ref{smooth Azumaya}.

It remains to prove that the maps on local rings $\O_{U,f(v)}\ra\O_{V,v}$ are bijective.
The schemes $U$ and $\GL_{n+2}$ are smooth over $\ZZ$, and  the same holds for the quotient
$V$. In particular,   $U$ is regular and $V$ is Cohen--Macaulay, such that
the finite morphism $f:V\ra U$ is flat, and we have to verify that it has degree one.
Since $U$ is connected, it suffices to show that the finite field extension $\kappa(f(v))\subset\kappa(v)$ has degree one for a single point $v\in V$.
This is indeed the case, because the morphism $f_\QQ$ is universally injective.
\qed

\begin{theorem}
\mylabel{clifford invariant generates}
Suppose that $n$ is odd. Let $S$ be a   integral noetherian scheme   whose strictly local rings
are factorial. Then 
the group $\Br(U_S)/\Br(S)$ is cyclic of order two, and the Clifford invariant   of the universal smooth quadric  
yields the generator.  
\end{theorem}

\proof
Let $F=\kappa(S)$ be the function field, and consider the commutative diagram
$$
\begin{CD}
0 @>>> \Br(S) @>>> \Br(U_S) @>>> \Br(U_S)/\Br(S) @>>> 0\\
@.     @VVV        @VVV          @VVV\\
0 @>>> \Br(F) @>>> \Br(U_F) @>>> \Br(U_F)/\Br(F) @>>> 0
\end{CD}
$$
whose rows are exact. The vertical maps on the left and middle are injective
by Proposition \ref{restriction function field}. Using the Snake Lemma and a diagram chase involving the compatiblesplitting of the
short exact sequences, we infer that the vertical map on the right is injective.
In turn, it suffices to treat the case $S=\Spec(k)$ the spectrum of a field. A similar
argument reduces to the case that $k$ is algebraically closed.
Now our task is to show that $\Br(U_k)$ is cyclic of order two, generated by the Clifford invariant
of the universal smooth quadric.

Set $G=\GL_{n+2,k}$ and $H=\GO_{n+2,k}$, which are smooth   connected   algebraic groups over $k$.
Throughout, we use the term \emph{algebraic group} for affine group scheme of finite type. 
The projection $G\ra G/H$ can be regarded as $H$-torsor.
We now use the identification of schemes $U_k=H\backslash G\simeq G/H$ from Proposition \ref{identification moduli}.
According to \cite{Sansuc 1981}, Proposition 6.10, there is an exact sequence
$$
\Pic(G)\stackrel{\varphi}{\lra}\Pic(H)\lra \Br(G/H)\lra\Br(G),
$$
where $\varphi$ is some homomorphism   defined in \cite{Sansuc 1981}, Lemma 6.4.
Using Proposition \ref{brauer gl} below, we conclude that $\Br(G/H)$ is cyclic of order two.

It remains to show that the Clifford invariant of the universal smooth  quadric is nonzero.
To this end, it suffices to exhibit a  single field extension $k\subset F$ together with 
a quadratic form $Q:F^{n+2}\ra F$ so that the resulting Azumaya algebra $A=\Cl_+(Q)$ 
has nontrivial Brauer class.  
This can be easily done with quaternion algebras. 
Let  $s,t$ be indeterminates, and consider the transcendental extension $F=k(s,t)$.
We denote by $q_0$ and $q_1$ the diagonal quadratic form 
given by $\langle s,t\rangle$ and $\langle1,1\rangle$, respectively. Furthermore, let $q_2=\langle 1\rangle$. 
Let $Q'$ be the direct sum
of $q_0$ and $r-1$ additional copies of $q_1$, where we write   $n+2=2r+1$. Set $Q=Q'\oplus q_2$.
Then
$$
\Cl_+(Q)=\Cl(Q') = \Cl(q_0)\otimes\Cl(q_1)\otimes\ldots\otimes\Cl(q_1).
$$
Each factor on the right can be viewed as a \emph{quaternion algebra}
$$
\left(\frac{a,b}{F}\right)=F\oplus Fi\oplus Fj\oplus Fk,\quad i^2=a,\; j^2=b,\; ij=k,
$$
namely $\Cl(q_0)=\left(\frac{s,t}{F}\right)$
and $\Cl(q_1)=\left(\frac{1,1}{F}\right)$. The former has nontrivial Brauer class, because
the \emph{Hilbert equation} $sx^2+ty^2=1$ has no solution in $F$, 
whereas the latter have trivial Brauer class (for example  \cite{Lam 2005}, Chapter III, Theorem 2.7).
\qed
 
\medskip
The result holds, in particular, if $S$ is the spectrum of a field $k$ or the ring of integers $\ZZ$.
From $\Br(\ZZ)=0$, it follows that the Brauer group $\Br(U)$ of the moduli space
of smooth quadrics over $\ZZ$ is cyclic of order two,   generated by the Clifford invariant
of the universal smooth quadric.
In the preceding proof, we have used the following facts:

\begin{proposition}
\mylabel{brauer gl}
Suppose $k$ is a an algebraically closed field. Let  $n\geq 1$ an   integer, $G=\GL_{n,k}$   the general linear group,
and $H=\GO_{n,k}$   the group of orthogonal similitudes.
Then $\Pic(G)=\Br(G)=0$.  Moreover, $\Pic(H)$ is cyclic of order two, provided
that $n$ is odd.
\end{proposition}

\proof
We may regard $G$ as an open subset inside $\AA^{d}_k$, $d=n^2$. The latter has
trivial Picard group and is locally factorial, whence $G$ has trivial Picard group.
We next verify the statement on the Brauer group.
Consider the special linear group $G'=\SL_{n,k}$, which sits in a short exact sequence
$
0\ra G'\ra G\stackrel{\det}{\ra}\GG_{m,k}\ra 0
$
of connected algebraic groups. This induces, by \cite{Sansuc 1981}, Corollary 6.11, an exact sequence
$$
\Br(G')\lra\Br(G)\lra\Br(\GG_{m,k}).
$$
The term on the right vanishes by Tsen's Theorem.
The term on the left vanishes as well, which can be seen as follows: There are no nontrivial characters 
$G'\ra\GG_{m,k}$. Moreover, the  algebraic fundamental group $\pi_1(G')$ 
 vanishes as well.
It thus follows from \cite{Iversen 1976}, Corollary 4.3 that $\Br(G')=0$.
Summing up,  $\Br(G)=0$.
Note that the \emph{algebraic fundamental group} of an algebraic group, which classifies central isogenies, usually differs
from the fundamental group of the underlying scheme, which classfies finite \'etale covering.

It remains to compute $\Pic(H)$ for $n$ odd. Let $H'=\SO_{n,k}$ be the special orthogonal group,
and consider the canonical morphism $H'\times\GG_{m,k}\ra H$.
This is an isomorphism, because $n$ is odd. Let  $\tilde{H}'=\Spin_{n,k}$ be the spin group,
and let $\tilde{H}'\ra H'$ be the canonical central isogeny of degree two,
whose kernel is isomorphic to the finite diagonalizable group scheme $D=\mu_{2,k}$.
According to \cite{Fossum; Iversen 1973}, Proposition 4.2, we have a short exact sequence
$$
\Hom(\tilde{H}',\GG_{m,k})\lra\Hom(D,\GG_{m,k}) \lra\Pic(H')\lra\Pic(\tilde{H}').
$$
The spin group $\tilde{H}'$ admits no nontrivial characters, whence the term on the left vanishes.
Obviously,  $\Hom(D,\GG_{m,k})$ is cyclic of order two. Moreover, the 
algebraic fundamental group $\pi_1(H')$ vanishes.
If follows that $\Pic(\tilde{H}')=0$, according to \cite{Fossum; Iversen 1973}, Corollary 4.5.
Summing up, $\Pic(H')$ is cyclic of order two.
Using the decomposition
$\Pic(H'\times\PP^1)=\Pic(H')\oplus\ZZ\O(1)$,
one easily infers   that the Picard group of $H=H'\times\GG_{m,k}\subset H'\times\PP^1$ is cyclic of order two.

\qed


\end{document}